\documentclass[11pt]{amsart}

\usepackage{amsmath}
\usepackage{amssymb}

\newtheorem{theorem}{Theorem}[section]
\newtheorem{claim}[theorem]{Claim}
\newtheorem{mclaim}[theorem]{Main Claim}

\newtheorem{lemma}[theorem]{Lemma}

\newtheorem{corollary}[theorem]{Corollary}

\theoremstyle{definition}
\newtheorem{definition}[theorem]{Definition}

\newtheorem{question}[theorem]{Question}

\theoremstyle{remark}
\newtheorem{remark}[theorem]{Remark}

\newcount\skewfactor
\def\mathunderaccent#1#2 {\let\theaccent#1\skewfactor#2
\mathpalette\putaccentunder}
\def\putaccentunder#1#2{\oalign{$#1#2$\crcr\hidewidth
\vbox to.2ex{\hbox{$#1\skew\skewfactor\theaccent{}$}\vss}\hidewidth}}
\def\name{\mathunderaccent\tilde-3 }

\def\smallbox#1{\leavevmode\thinspace\hbox{\vrule\vtop{\vbox
   {\hrule\kern1pt\hbox{\vphantom{\tt/}\thinspace{\tt#1}\thinspace}}
   \kern1pt\hrule}\vrule}\thinspace}

\newcommand{\cf}{{\rm cf}}

\def\qedref#1{$\qed_{\reforiginal{#1}}$}

\title{Martin's maximum and the non-stationary ideal}

\author{Shimon Garti}
\address{Institute of Mathematics,
 The Hebrew University of Jerusalem,
 Jerusalem 91904, Israel}
\email{shimon.garty@mail.huji.ac.il}

\author{Yair Hayut}
\address{Institute of Mathematics,
 The Hebrew University of Jerusalem,
 Jerusalem 91904, Israel}
\email{yair.hayut@mail.huji.ac.il}

\author{Haim Horowitz}
\address{Institute of Mathematics,
 The Hebrew University of Jerusalem,
 Jerusalem 91904, Israel}
\email{haim.horowitz@mail.huji.ac.il}

\author{Menachem Magidor}
\address{Institute of Mathematics,
 The Hebrew University of Jerusalem,
 Jerusalem 91904, Israel}
\email{menachem@math.huji.ac.il}

\thanks{}


\begin{document}
\let\labeloriginal\label
\let\reforiginal\ref

\begin{abstract}
We deal with ${\rm NS}_{\omega_1}$ and $\mathcal{D}_{\omega_1}$ under Martin's maximum and the proper forcing axiom.
\end{abstract}

\maketitle

\newpage

\section{Introduction}

In this paper we try to study the influence of forcing axioms on the non-stationarty ideal over $\aleph_1$ (denoted by ${\rm NS}_{\omega_1}$) and its dual, the club filter $\mathcal{D}_{\omega_1}$.
In the first section we introduce a general reflection theorem from which we derive a saturation property of ${\rm NS}_{\omega_1}$.
In the second section we deal with Galvin's property with respect to $\mathcal{D}_{\omega_1}$.

\begin{definition}
\label{defsat} Saturation. \newline
An ideal $\mathcal{J}$ over $\kappa$ is $(\mu,\lambda,\theta)$-saturated iff for every collection $\{A_\alpha:\alpha<\mu\}\subseteq\mathcal{J}^+$ there is a sub-collection $\{A_{\alpha_\varepsilon}:\varepsilon\in\lambda\}$ such that $\bigcap\{A_{\alpha_\varepsilon}:\varepsilon\in y\}\in\mathcal{J}^+$ for every $y\in[\lambda]^\theta$.
\end{definition}

An ideal $\mathcal{I}$ over a successor cardinal $\kappa^+$ is Laver iff $\mathcal{I}$ is uniform (i.e., contains all the bounded subsets of $\kappa^+$), $\kappa^+$-complete and $(\kappa^{++},\kappa^{++},\kappa)$-saturated. If we require only $(\kappa^{++},\kappa^+,\kappa)$-saturation then $\mathcal{I}$ is called weakly Laver. If we have just $(\kappa^{++},\kappa^+,<\kappa)$-saturation then $\mathcal{I}$ is \emph{almost weakly Laver}.

Laver ideals can be forced over $\aleph_1$ by collapsing sufficiently large cardinals. In the original work of Laver, \cite{MR673792}, a huge cardinal was used. Other concepts of saturation can be forced with the same basic idea, i.e. to collapse a large cardinal and preserve some of its qualities.
An interesting question is what happens under classical forcing axioms. One expects $\aleph_1$ to behave like a large cardinal under these axioms.

The simplest notion of saturation for an ideal $\mathcal{I}$ over $\kappa^+$ says that if $\{A_\alpha:\alpha<\kappa^{++}\}\subseteq\mathcal{I}^+$ then there are $\alpha<\beta<\kappa^{++}$ so that $A_\alpha\cap A_\beta\in\mathcal{I}^+$. Such an ideal is called Kunen, and the following theorem from \cite{MR924672} embodies the above expectation:

\begin{theorem}
\label{thmkunen} Under MM, the ideal ${\rm NS}_{\omega_1}$ is Kunen.
\end{theorem}

\hfill \qedref{thmkunen}

So strong forcing axioms yield saturation properties at $\aleph_1$. However, the situation is more involved. Laverness, for example, is excluded from $\aleph_1$ under any classical forcing axiom. The theorem below is typical, its first assertion is due to Larson, \cite{MR2146222}, and its second part appears in \cite{Ga}.

\begin{theorem}
\label{thmmaba} Let $\mathcal{J}$ be a uniform ideal.
\begin{enumerate}
\item [$(\aleph)$] If $MA + 2^\omega=\omega_2$ and $\mathcal{J}$ is an ideal over $\aleph_1$ then $\mathcal{J}$ is not Laver.
\item [$(\beth)$] If $BA + 2^{\omega_1}=\omega_3$ and $\mathcal{J}$ is an ideal over $\aleph_2$ then $\mathcal{J}$ is not Laver.
\end{enumerate}
\end{theorem}

\hfill \qedref{thmmaba}

A natural question, raised by Larson, is about the intermediate concept of weak Laverness at $\aleph_1$ under some forcing axiom. Larson asked, in particular, about ${\rm NS}_{\omega_1}$.
In the first section we obtain the following result:

\begin{theorem}
\label{thmawlaver} If MM holds then ${\rm NS}_{\omega_1}$ is almost weakly Laver.
\end{theorem}

\hfill \qedref{thmawlaver}

We indicate that this result follows in a different way from \cite{MR3744886}, Lemma 3.10.
In the second section we address the cardinal characteristic $\mathfrak{gp}$, related to $\mathcal{D}_{\omega_1}$.
We also deal with Galvin's property under the proper forcing axiom.
Regarding this property we prove the consistency of PFA with the strong negation of Galvin's property.
We believe that this strong negation holds in any model of PFA.

We shall use the standard notation. If $\mathcal{I}$ is an ideal over $\kappa^+$ then $\mathcal{I}$ will be uniform and $\kappa^+$-complete unless otherwise stated. The set $\mathcal{P}(\kappa^+)\setminus\mathcal{I}$ is denoted by $\mathcal{I}^+$ (the $\mathcal{I}$-positive sets). We employ the Jerusalem forcing notation, i.e. $p\leq q$ reads $p$ is weaker than $q$.

If $g:\kappa\rightarrow S$ and $\beta\in\kappa$ then $g''\beta = \{g(\alpha):\alpha<\beta\}$.
The set $\mathcal{P}_\kappa\lambda$ is the collection of subsets of $\lambda$ whose cardinality is less than $\kappa$.
Let $M,N$ be models of (enough) ZFC. We say that $N$ is a $\kappa$-end-extension of $M$ iff for every $\zeta\in N\cap\kappa$ either $\zeta\in M$ or $\zeta\geq \sup(M\cap\kappa)$. This concept will be used below with respect to $\kappa=\aleph_2$.

Let $A$ be a set and let $F:[A]^{<\omega}\rightarrow\mathcal{P}_\kappa(A)$ be a function. An element $x$ of $\mathcal{P}_\kappa(A)$ will be called a closure point of $F$ iff $F(e)\subseteq x$ for every $e\in[x]^{<\omega}$. The set of closure points of $F$ will be denoted by $c\ell_F$. It is straightforward that this set is a club in $\mathcal{P}_\kappa(A)$. Moreover, every club $C$ of $\mathcal{P}_\kappa(A)$ can be coded by a function $F$ in the sense that there exists a function $F:[A]^{<\omega}\rightarrow\mathcal{P}_\kappa(A)$ such that $c\ell_F\subseteq C$.

A central concept from \cite{MR924672} is $\aleph_1$-reflection. Suppose that $S\subseteq[\mathcal{H}(\theta)]^\omega$. We shall say that $S$ reflects to a set of size $\aleph_1$ iff there exists $x\subseteq \mathcal{H}(\theta), |x|=\aleph_1, \omega_1\subseteq x$ such that $S\cap[x]^\omega$ is stationary in $[x]^\omega$. The following is labeled in \cite{MR924672} as Theorem 13:

\begin{theorem}
\label{thmreflection} MM and $\aleph_1$-reflection. \newline
Assume MM and let $\theta$ be a regular uncountable cardinal. \newline
Every stationary set $S\subseteq [\mathcal{H}(\theta)]^\omega$ reflects to a set of size $\aleph_1$.
\end{theorem}

\hfill \qedref{thmreflection}

\newpage

\section{Almost weak Laverness}

Our first statement is a strengthening of typical properties which can be proved under MM. It says, roughly, that we have a good control on $M\cap\omega_1$ for countable elementary submodels of $\mathcal{H}(\theta)$. We shall need, however, a simple lemma ahead of the main claim.

\begin{lemma}
\label{lemchain} Let $\mathbb{Q}$ be an $\aleph_2$-cc forcing notion, and let $G\subseteq\mathbb{Q}$ be generic.
Let $\{q_\varepsilon:\varepsilon\in\omega_2\}$ be any collection of conditions in $\mathbb{Q}$. \newline
There exists a condition $p\in\mathbb{Q}$ which forces $\aleph_2$-many members from this collection into the generic set $G$.
\end{lemma}

\par\noindent\emph{Proof}. \newline
Fix a collection $\{q_\varepsilon:\varepsilon\in\omega_2\}\subseteq \mathbb{Q}$.
We shall prove that some $q_\varepsilon$ forces that $|\{\zeta:q_\zeta \in \name{G}\}|=\aleph_2$.
Toward contradiction assume that there is no such condition.
Hence for every $\varepsilon\in\omega_2$ there are $p_\varepsilon\in\mathbb{Q}$ and an ordinal $\alpha_\varepsilon$ such that $q_\varepsilon\leq p_\varepsilon$ and $p_\varepsilon\Vdash \alpha_\varepsilon$ is the first ordinal such that $q_\zeta\notin\name{G}$ for every $\zeta\in[\alpha_\varepsilon,\omega_2)$. Let $\delta = \bigcup_{\varepsilon\in\omega_2}\alpha_\varepsilon$.

\par\noindent\emph{Case 1}: $\delta=\omega_2$.

In this case there is an increasing sequence $\langle \varepsilon_\zeta: \zeta\in\omega_2\rangle$ of ordinals in $\omega_2$ such that $\langle \alpha_{\varepsilon_\zeta}:\zeta\in\omega_2\rangle$ is strictly increasing. We claim that the corresponding set $\{p_{\alpha_{\varepsilon_\zeta}} : \zeta\in\omega_2\}$ is an antichain in $\mathbb{Q}$, contradicting the $\aleph_2$-chain condition.
For proving this fact assume that $\eta<\zeta<\omega_2$. So $p_{\alpha_{\varepsilon_\eta}}\Vdash q_\zeta\notin\name{G}$ for every $\zeta\in[\alpha_{\varepsilon_\eta},\omega_2)$ while $p_{\alpha_{\varepsilon_\zeta}}\Vdash \exists\tau\in [\alpha_{\varepsilon_\eta},\omega_2), q_\tau\in\name{G}$. It means that $p_{\alpha_{\varepsilon_\eta}}\perp p_{\alpha_{\varepsilon_\zeta}}$, so we are done.

\par\noindent\emph{Case 2}: $\delta<\omega_2$.

In this case there is an ordinal $\alpha<\omega_2$ such that $|\{\varepsilon\in\omega_2: \alpha_\varepsilon=\alpha\}|=\aleph_2$. Choose any $\zeta\in\omega_2$ such that $\alpha<\zeta$ and $\alpha_\zeta=\alpha$.
On the one hand $p_\zeta\Vdash q_\zeta\in\name{G}$ since $q_\zeta\leq p_\zeta$. On the other hand, $p_\zeta$ forces that $q_\tau\notin\name{G}$ whenever $\tau\in[\alpha_\zeta,\omega_2)=[\alpha,\omega_2)$. In particular, $p_\zeta \Vdash q_\zeta\notin\name{G}$, a contradiction.

\hfill \qedref{lemchain}

The main claim, to be proved below, says that one can find \emph{nice extensions} of countable elementary submodels. We shall prove it in two steps, first relativized to some club subset and then in a general form. Ahead of the claims (and proofs) we define the properties that we shall need.

\begin{definition}
\label{defnice} Nice extensions. \newline
Assume that $M$ is countable and $M\prec\mathcal{H}(\theta)$ for some regular cardinal $\theta$.
Assume further that $\{A_\beta:\beta\in\omega_2\}\subseteq {\rm NS}^+_{\omega_1}$.
An elementary submodel $N\prec\mathcal{H}(\theta)$ will be called a nice extension of $M$ iff:
\begin{enumerate}
\item [$(\aleph)$] $M\subseteq N$.
\item [$(\beth)$] $N$ is countable.
\item [$(\gimel)$] $N$ is an $\omega_2$-end-extension of $M$.
\item [$(\daleth)$] $\exists\beta\in N\setminus M$ such that $M\cap\omega_1 = N\cap\omega_1\in A_\beta$.
\end{enumerate}
\end{definition}

The third requirement, namely $N$ is an $\omega_2$-end-extension of $M$, can be verified while concentrating on ordinals of $\omega_1$. This is the content of the following:

\begin{lemma}
\label{lemextension} Suppose that $\theta=\cf(\theta)$ and $M,N$ are countable structures which satisfy $M\prec N\prec \mathcal{H}(\theta)$. \newline
If $M\cap\omega_1=N\cap\omega_1$ then $N$ is an $\omega_2$-end-extension of $M$.
\end{lemma}

\par\noindent\emph{Proof}. \newline
Assume that $\xi\in N\cap\omega_2$ and $\zeta\in M\cap\omega_2$ for some $\zeta>\xi$. We need showing that $\xi\in M$. If $\xi\in\omega_1$ then $\xi\in M$ as we are assuming that $M\cap\omega_1 = N\cap\omega_1$. So we may assume that $\omega_1\leq\xi$, and hence $\zeta>\omega_1$.

Now $\zeta\in\omega_2\cap M$ and $\zeta>\omega_1$ so there is a bijection $h:\zeta\rightarrow\omega_1, h\in M$. Since $\xi<\zeta$ we see that $h(\xi)\in\omega_1$. Since $\xi\in N$ and $h\in M\subseteq N$ we see that $h(\xi)\in N$. Therefore, $h(\xi)\in N\cap\omega_1 = M\cap\omega_1$, hence $h(\xi)\in M$ as well. By elementarity $h^{-1}\in M$ and hence $\xi=h^{-1}(h(\xi))\in M$, so we are done.

\hfill \qedref{lemextension}

We are ready now for the first assertion:

\begin{claim}
\label{clmclub} Assume MM. \newline
Suppose that $\{A_\gamma: \gamma\in\omega_2\}\subseteq {\rm NS}^+_{\omega_1}, A\in {\rm NS}^+_{\omega_1}$ and every stationary subset of $A$ intersects $\omega_2$-many $A_\gamma$-s positively. \newline
Then for every sufficiently large regular cardinal $\theta$ there is a club $C\subseteq \mathcal{P}_{\omega_1}(\mathcal{H}(\theta))$ which consists of countable elementary submodels of $\mathcal{H}(\theta)$ such that for every $M\in C$ if $M\cap\omega_1\in A$ then there exists a nice extension $N$ of $M$.
\end{claim}

\par\noindent\emph{Proof}. \newline
Toward contradiction assume that such a club $C$ does not exist. Fix a stationary subset $S$ of $\mathcal{P}_{\omega_1}(\mathcal{H}(\theta))$ such that every $M\in S$ is a countable elementary submodel of $\mathcal{H}(\theta)$ satisfying $M\cap\omega_1\in A$ and there is no nice extension of $M$.

By the reflection principle of Theorem \ref{thmreflection} there is some $I\subseteq\mathcal{H}(\theta), |I|=\aleph_1, \omega_1\subseteq I$ such that $S\cap\mathcal{P}_{\omega_1}(I)$ is stationary in $\mathcal{P}_{\omega_1}(I)$.
Fix a bijection $g:\omega_1\rightarrow I$. Define $T = \{\beta\in\omega_1: g''\beta\in S\}$. We claim that $T$ is a stationary subset of $\omega_1$, and moreover $T\cap A$ is stationary. For suppose not, and choose a club $D\subseteq \omega_1$ such that $D\cap(T\cap A)=\varnothing$. Without loss of generality if $\alpha\in D$ then $\omega_1\cap g''\alpha=\alpha$, as we can cut $D$ with the club of ordinals $\alpha$ for which $\omega_1\cap g''\alpha=\alpha$.

Let $E = \{M\in \mathcal{P}_{\omega_1}(I): M\cap\omega_1\in D, g''(M\cap\omega_1)\subseteq M \wedge g^{-1}[M]\subseteq M\}$. Notice that $E$ is a club in $\mathcal{P}_{\omega_1}(I)$. Likewise, $M=g''(M\cap\omega_1)$ whenever $M\in E$. For this, if $x\in M$ then $x=g(g^{-1}(x))$. But $g^{-1}(x)\in M\cap\omega_1$ since $M$ is closed under $g^{-1}$, so $g(g^{-1}(x))\in g''(M\cap\omega_1)$. It means that $M\subseteq g''(M\cap\omega_1)$, and $g''(M\cap\omega_1)\subseteq M$ follows from the fact that $M$ is closed under $g$.

Pick any element $M\in E\cap S$. On the one hand $M\cap\omega_1\in D$ since $M\in E$. On the other hand $M\cap\omega_1\in T\cap A$ since $M\in S$ and $M=g''(M\cap\omega_1)$. In other words, the ordinal $M\cap\omega_1$ belongs to $D\cap(T\cap A)$, a contradiction.

Define $\eta=\sup(I\cap\omega_2)$. Let $D' = \{\beta\in\omega_1: \omega_1\cap g''\beta = \beta\}$, and let $B' = (T\cap A)\cap D'$. The set $B'\subseteq A$ is stationary, and by the properties of $A$ there is an ordinal $\gamma\in\omega_2$ such that $\eta<\gamma$ and $B = B'\cap A_\gamma$ is stationary in $\omega_1$. Choose $N\prec\mathcal{H}(\theta), |N|=\aleph_0$ such that $g,\gamma\in N$ and $\beta=N\cap\omega_1\in B$. Set $M_0=g''\beta$. Notice that $M_0\in S$ since $\beta\in B\subseteq T$. We claim that $N$ is a nice extension of $M_0$, which is an absurd as $M_0\in S$.

First observe that $M_0\subseteq N$ since $g\in N$ and $\beta\subseteq N$. Second, $N\cap\omega_1=\beta=M_0\cap\omega_1$ since $\beta\in D'$, so Lemma \ref{lemextension} ensures that $N$ is an $\omega_2$-end-extension of $M_0$. Finally, $\gamma\in N$ by the choice of $N$ but $\gamma\notin M_0$. Indeed, $M_0=g''\beta\subseteq I$ so if $\gamma\in M_0$ then $\gamma=g(\alpha)$ for some $\alpha<\beta$. It means that $\gamma\in I$ (as $g:\omega_1\rightarrow I$) and $\gamma\in\omega_2$, which is impossible since $\gamma>\eta=\sup(I\cap\omega_2)$.
Summing up, $N$ is a nice extension of the element $M_0\in S$, a contradiction.

\hfill \qedref{clmclub}

We proceed to the second step, which is basically the same statement but without applying to the club $C\subseteq \mathcal{P}_{\omega_1}(\mathcal{H}(\theta))$.

\begin{mclaim}
\label{mclm} Nice extensions. \newline
Assume MM and Suppose that $\{A_\gamma: \gamma\in\omega_2\}\subseteq {\rm NS}^+_{\omega_1}, A\in {\rm NS}^+_{\omega_1}$, and for every $S\in {\rm NS}^+_{\omega_1}, S\subseteq A$ there are $\omega_2$-many ordinals $\beta$ for which $S\cap A_\beta\in {\rm NS}^+_{\omega_1}$. \newline
Then for every large enough regular cardinal $\chi$ and every countable $M\prec\mathcal{H}(\chi)$ so that $M\cap\omega_1\in A$ one can find a nice extension $N\prec\mathcal{H}(\chi)$.
\end{mclaim}

\par\noindent\emph{Proof}. \newline
We shall assume that $\mathcal{H}(\chi)$ is augmented by a fixed well-ordering and any required functions or predicates, e.g. definable Skolem functions and the sequence $\{A_\gamma: \gamma\in\omega_2\}$ as well as the set $A$ as predicates. Let $\theta=\cf(\theta)$ be large enough to satisfy Claim \ref{clmclub}. We may assume that every function from the finite subsets of $\omega_2$ into $\omega_2$ belongs to $\mathcal{H}(\theta)$, upon noticing that this requires $\theta\geq\omega_3$. Choose $\chi=\cf(\chi)>\theta$ so that $\mathcal{H}(\theta)\in \mathcal{H}(\chi)$. The latter is needed in order to guarantee that the club $C$ from Claim \ref{clmclub} belongs to $\mathcal{H}(\chi)$.

Assume that $M\prec\mathcal{H}(\chi), |M|=\aleph_0$ and $M\cap\omega_1\in A$. It follows that $M\cap\mathcal{H}(\theta)\in \mathcal{P}_{\omega_1}(\mathcal{H}(\theta))$, and we may assume that $M\cap \mathcal{H}(\theta)\in C$. Indeed, we can assume that the club $C$ is coded by a function $F$ which appears as a predicate in the augmented structure of $\mathcal{H}(\chi)$. Now $F$ is a function from $\mathcal{H}(\theta)^{<\omega}$ into $\mathcal{H}(\theta)$ and $M$ is an elementary submodel of $\mathcal{H}(\chi)$ with this predicate, so together $M\cap\mathcal{H}(\theta)$ is closed under $F$ and hence belongs to $C$.

By virtue of Claim \ref{clmclub} there is a nice extension $N_0$ of $M\cap \mathcal{H}(\theta)$, and in particular $N_0\prec \mathcal{H}(\theta)$. We wish to extend $N_0$ to an elementary submodel of $\mathcal{H}(\chi)$ which will be a nice extension of $M$.

Define $N={\rm Sk}^{\mathcal{H}(\chi)}(M\cup(N_0\cap\omega_2))$. Notice that $N\prec\mathcal{H}(\chi)$, and we shall prove that $N$ is as desired. Being the Skolem hull of a countable set, $N$ is virtually countable. Clearly, $M\subseteq N$. Likewise, there exists an ordinal $\beta\in N_0\setminus(M\cap \mathcal{H}(\theta))$ such that $(M\cap \mathcal{H}(\theta))\cap\omega_1 = N_0\cap\omega_1\in A_\beta$. Since $\theta$ is large enough, $(M\cap \mathcal{H}(\theta))\cap\omega_1 = M\cap\omega_1$, and $\beta\in N_0\cap\omega_2\subseteq N$. So if we prove that $N_0\cap\omega_2 = N\cap\omega_2$ we will conclude that $N_0\cap\omega_1 = N\cap\omega_1\in A_\beta$. Moreover, $N_0\cap\omega_2 = N\cap\omega_2$ accomplishes all the requirements listed in the definition of a nice extension.

Let us prove, therefore, that $N_0\cap\omega_2 = N\cap\omega_2$. Clearly, $N_0\cap\omega_2 \subseteq N\cap\omega_2$, so assume that $\xi\in N\cap\omega_2$ with the goal of showing that $\xi\in N_0$ as well.
By the definition of $N$ there exists a Skolem function $f$, a finite set $\{a_1,\ldots,a_k\}\subseteq M$ and a finite set of ordinals $\{\gamma_1,\ldots,\gamma_\ell\}\subseteq N_0\cap\omega_2$ such that:
$$
\xi = f(a_1,\ldots,a_k,\gamma_1,\ldots,\gamma_\ell).
$$
We may assume that ${\rm Rang}(f)\subseteq\omega_2$ (if not, replace $f$ by $g$ which agrees with $f$ whenever $f(e)\in\omega_2$ and assumes zero otherwise). Define $h:[\omega_2]^\ell\rightarrow\omega_2$ by $h(t) = f(a_1,\ldots,a_k,t)$.

We claim that $h\in M$.
Indeed, $M\prec\mathcal{H}(\chi)$ and $\mathcal{H}(\chi)$ contains a complete set of Skolem functions, including the specific function $f$. The function $h$ is definable from $f$ and $\{a_1\ldots,a_k\}\subseteq M$ hence belongs to $M$.
Likewise, assuming that $\theta$ is sufficiently large we may assume that all the functions from $[\omega_2]^\ell$ into $\omega_2$ are in $\mathcal{H}(\theta)$, in particular $h\in \mathcal{H}(\theta)$ and hence $h\in M\cap \mathcal{H}(\theta)\subseteq N_0$. Since $\{\gamma_1,\ldots,\gamma_\ell\}\subseteq N_0$ we conclude that $\xi = f(a_1,\ldots,a_k,\gamma_1,\ldots,\gamma_\ell) = h(\gamma_1,\ldots,\gamma_\ell)\in N_0$, and the proof is accomplished.

\hfill \qedref{mclm}

\begin{remark}
\label{rmeasurable} The ability to produce nice extensions of countable elementary submodels resembles similar theorems which hold at measurable cardinals. It means that $\aleph_2$ behaves like a measurable cardinal under MM. This is coherent with other statements which show that small accessible cardinals possess large cardinal properties under strong forcing axioms.
\end{remark}

\hfill \qedref{rmeasurable}

As indicated in the introduction, the fact that ${\rm NS}_{\omega_1}$ is almost weakly Laver under MM follows from Lemma 3.10 of \cite{MR3744886}.
It can also be derived from the main claim above:

\begin{theorem}
\label{thmmt} The ideal ${\rm NS}_{\omega_1}$ is $(\aleph_2,\aleph_1,<\aleph_0)$-saturated under MM.
\end{theorem}

\par\noindent\emph{Proof}. \newline
Suppose that $\{A_\alpha:\alpha\in\omega_2\}\subseteq {\rm NS}^+_{\omega_1}$. Let $A\subseteq\omega_1$ be as guaranteed by Lemma \ref{lemchain}. By induction on $\gamma\in\omega_1$ we choose (for sufficiently large $\chi$) a model $N_\gamma\prec\mathcal{H}(\chi)$ as follows:

$\gamma=0$: Let $N_0$ be a countable elementary submodel of $\mathcal{H}(\chi)$ so that $N_0\cap\omega_1\in A$. We also assume that $\{A_\alpha:\alpha\in\omega_2\}\in N_0$.

$\gamma+1$: Denote $N_\gamma$ by $M$, and choose $N$ as described in Claim \ref{mclm} with respect to $M$. Let $\alpha_\gamma\in N\setminus M$ be the ordinal given by the definition of nice extensions. Set $N_{\gamma+1} = N$.

$\gamma$ is limit: Let $N_\gamma = \bigcup_{\delta<\gamma}N_\delta$.

Denote the union $\bigcup_{\gamma\in\omega_1}N_\gamma$ by $N$. Notice that $\{A_{\alpha_\gamma}:\gamma\in\omega_1\}\subseteq N$. We claim that this collection exemplifies the almost weak Laverness of ${\rm NS}_{\omega_1}$. For this, let $u$ be a finite subset of $\omega_1$ and $S = \bigcap_{\gamma\in u}A_{\alpha_\gamma}$. Notice that $S\in N$ as $u$ is finite.

Assume toward contradiction that $S\in{\rm NS}_{\omega_1}$, and let $C$ be a club subset of $\omega_1$ so that $C\in N$ and $N\models C\cap S=\varnothing$. Let $\tau$ be the characteristic ordinal $N_0\cap\omega_1$. By the nature of the $N_\gamma$-s we see that $N_\gamma\cap\omega_1=\tau$ for every $\gamma\in\omega_1$ and hence $N\cap\omega_1=\tau$.

On the one hand, $\tau\in C$ since $\mathcal{H}(\chi)$ believes that $C$ is unbounded below $\tau$ and $C$ is closed. On the other hand, $\tau\in A_{\alpha_\gamma}$ for each $\gamma\in u$ since $N_\gamma\cap\omega_1 = N\cap\omega_1 = \tau$ for every $\gamma\in\omega_1$ and by the definition of nice extensions. In particular $\tau\in S$, a contradiction.

\hfill \qedref{thmmt}

We make the comment that the finite number of $A_{\alpha_\beta}$-s is required only in order to make sure that $S = \bigcap_{\beta\in u}A_{\alpha_\beta}\in N$. This seems to be an essential obstacle while considering the possibility of strengthening the above theorem. Indeed, no closure assumption on $N$ can be made (apart from the closure to finite sequences which follows from elementarity).

The question of weak Laverness under MM remains open, but almost weak Laverness is quite settled. It is interesting to compare MM with PFA at this point:

\begin{corollary}
\label{corawl} Forcing axioms and the non-stationary ideal.
\begin{enumerate}
\item [$(\aleph)$] ${\rm NS}_{\omega_1}$ is almost weakly Laver under MM.
\item [$(\beth)$] The almost weak Laverness of ${\rm NS}_{\omega_1}$ is independent over the PFA.
\end{enumerate}
\end{corollary}

\par\noindent\emph{Proof}. \newline
The first statement is the above theorem, and PFA is consistent with ${\rm NS}_{\omega_1}$ being almost weakly Laver since MM implies PFA. The consistency of PFA with the failure of almost weak Laverness of ${\rm NS}_{\omega_1}$ and even with the failure of the simplest saturation property is due to Velickovic, \cite{vetemp} (and as noted there, this result was obtained independently by Shelah).

\hfill \qedref{corawl}

\begin{remark}
\label{r0} For the main result we used MM, but actually almost weak Laverness follows from the $\aleph_1$-reflection principle and the fact that ${\rm NS}_{\omega_1}$ is Kunen. This might be meaningful when other axioms are deemed with respect to various saturation properties.
\end{remark}

\newpage

\section{The club filter}

A celebrated theorem of Galvin which appeared in \cite{MR0369081} says that under the continuum hypothesis any collection of $\aleph_2$ many clubs of $\aleph_1$ contains a sub-collection of size $\aleph_1$ whose intersection is a club subset of $\aleph_1$.
Remark that the collection of end-segments of $\aleph_1$ is a family of $\aleph_1$ many clubs of $\aleph_1$, each sub-family of which of size $\aleph_1$ has an empty intersection.
This means that one can gather a small number of clubs of $\aleph_1$ with empty intersection of each $\aleph_1$ of them, but (under the continuum hypothesis) one cannot collect many clubs in this way. Here, many clubs means $\aleph_2$ clubs of $\aleph_1$.

Galvin's result can be generalized in several directions, and one of them is based on replacing the continuum hypothesis by Devlin-Shelah's weak diamond.
It has been proved in \cite{MR3604115} that if $2^\omega<2^{\omega_1}$ (this is equivalent to the weak diamond at $\aleph_1$) then any collection of $\lambda$ clubs of $\aleph_1$ contains a sub-collection of size $\aleph_1$ whose intersection is a club subset of $\aleph_1$ where $\lambda=(2^\omega)^+$.
This is quite a surprising result, since $2^{\omega_1}$ can be much larger than $\lambda$ and yet Galvin's property holds here with respect to any collection of $\lambda$ sets. Actually, this gives rise to defining the first cardinal which satisfies this property as a cardinal characteristic of the continuum:

\begin{definition}
\label{defgp} The cardinal characteristic $\mathfrak{gp}$. \newline
We define $\mathfrak{gp}$ as the minimal $\kappa$ such that every family $\{C_\alpha:\alpha\in\kappa^+\}$ of club subsets of $\aleph_1$ contains a subfamily $\{C_{\alpha_\beta}:\beta\in\omega_1\}$ whose intersection is closed and unbounded in $\aleph_1$.
\end{definition}

The definition comes from \cite{MR3787522} and it has been shown there that $\mathfrak{gp}$ is well-defined and assumes a value between $\aleph_1$ and $2^{\aleph_0}$.
The behavior of $\mathfrak{gp}$ when related to other cardinal characteristics is very peculiar.
It is consistent that $\mathfrak{gp}<\mathfrak{m}$, but it is also consistent that $\mathfrak{gp}>\mathfrak{i}$.
Dealing with $\mathfrak{gp}$ by itself, it has been asked in \cite{MR3787522} what are the possible cofinalities of $\mathfrak{gp}$ and, in particular, is it consistent that $\cf(\mathfrak{gp})=\omega$.
In this paper we try to deal with this question.
Like many other cardinal characteristics, the definition of $\mathfrak{gp}$ generalizes to higher cardinals.
This is central in our context, since our knowledge on $\cf(\mathfrak{gp}_\kappa)$ for $\kappa=\cf(\kappa)\geq\aleph_3$ is a bit better.

\begin{definition}
\label{defgpkappa} The cardinal characteristic $\mathfrak{gp}_\kappa$. \newline
Let $\kappa=\cf(\kappa)$. \newline
We define $\mathfrak{gp}_\kappa$ as the minimal $\lambda$ such that every family $\{C_\alpha:\alpha\in\lambda^+\}$ of club subsets of $\kappa^+$ contains a subfamily $\{C_{\alpha_\beta}:\beta\in\kappa^+\}$ whose intersection is closed and unbounded in $\kappa^+$.
\end{definition}

A straightforward modification of the statements in \cite{MR3787522} shows that $\mathfrak{gp}_\kappa$ is well-defined and that $\kappa^+\leq\mathfrak{gp}_\kappa\leq 2^\kappa$.
For basic background regarding pcf theory we suggest \cite{MR2768693}, and for further information we refer to \cite{MR1318912}.

We commence with showing that under some pcf assumptions one can prove that $\cf(\mathfrak{gp})>\omega$. The assumption is a consequence of the existence of a good scale, which in turn can be proved to exist from instances of weak square. However, this assumption cannot be established in ZFC, and actually it fails if one assumes certain instances of Chang's conjecture.
We shall prove our first theorem based on a pcf assumption, and later we shall connect this assumption to known combinatorial concepts like weak squares and Chang's conjecture.

\begin{theorem}
\label{thmcofgp}
Let $\mu$ be a singular cardinal.
\begin{enumerate}
\item [$(A)$] Suppose that:
\begin{enumerate}
\item [$(a)$] $\omega=\cf(\mu)<\mu<2^\omega$ and $(\mu_n:n\in\omega)$ is an increasing sequence of regular cardinals such that $\mu=\bigcup_{n\in\omega}\mu_n$.
\item [$(b)$] $J\supseteq J^{\rm bd}_\omega$ and ${\rm tcf}(\prod_{n\in\omega}\mu_n,J)=\mu^+$, as witenessed by the scale $(f_\alpha:\alpha\in\mu^+)$.
\item [$(c)$] For every $A\subseteq\mu^+$ such that $|A|=\aleph_1$ it is true that $|\{f_\alpha(n):\alpha\in A,n\in\omega\}|=\aleph_1$.
\end{enumerate}
Then $\mathfrak{gp}\neq\mu$.
\item [$(B)$] Suppose that:
\begin{enumerate}
\item [$(a)$] $\omega_1=\cf(\mu)<\mu\leq 2^\omega$ and $(\mu_\gamma:\gamma\in\omega_1)$ is an increasing sequence of regular cardinals such that $\mu=\bigcup_{\gamma\in\omega_1}\mu_\gamma$.
\item [$(b)$] $J\supseteq J^{\rm bd}_{\omega_1}$ and ${\rm tcf}(\prod_{\gamma\in\omega_1}\mu_\gamma,J)=\mu^+$, as witenessed by the scale $(g_\alpha:\alpha\in\mu^+)$.
\item [$(c)$] For every $A\subseteq\mu^+$ such that $|A|=\aleph_1$ there is $\delta\in\omega_1$ such that $|\{g_\alpha(\gamma):\alpha\in A,\gamma\in\delta\}|=\aleph_1$.
\end{enumerate}
Then $\mathfrak{gp}\neq\mu$.
\end{enumerate}
\end{theorem}

\par\noindent\emph{Proof}. \newline
For the first part, assume toward contradiction that $\mathfrak{gp}=\mu$.
Fix a collection $\mathcal{D}=\{D_\alpha:\alpha\in\mu\}$ such that $\bigcap\{D_{\alpha_\beta}:\beta\in B\}$ is bounded in $\omega_1$ for every $B\in[\mu]^{\aleph_1}$.
For every $\alpha\in\mu^+$ let $C_\alpha=\bigcap\{D_{f_\alpha(n)}:n\in\omega\}$.
Each $C_\alpha$ is a club subset of $\aleph_1$, being the intersection of but $\aleph_0$-many clubs of $\aleph_1$.
Set $\mathcal{C}=\{C_\alpha:\alpha\in\mu^+\}$.

Since $\mathfrak{gp}=\mu$, one can choose $A\subseteq\mu^+,|A|=\aleph_1$ such that $C=\bigcap\{C_\alpha:\alpha\in A\}$ is a club subset of $\aleph_1$.
By our assumptions, $|\{f_\alpha(n):\alpha\in A, n\in\omega\}|=\aleph_1$ and hence $\{D_{f_\alpha(n)}:\alpha\in A,n\in\omega\}\subseteq [\mathcal{D}]^{\aleph_1}$.
By the nature of the collection $\mathcal{D}$, the set $D=\bigcap\{D_{f_\alpha(n)}:\alpha\in A, n\in\omega\}$ is bounded in $\omega_1$.
However, $D=C$, a contradiction.

For the second part, assume toward contradiction that $\mathfrak{gp}=\mu$ and fix a family of clubs $\mathcal{D}=\{D_\alpha:\alpha\in\mu\}$ which exemplifies this fact.
For every $\alpha\in\mu^+$ we define $C_\alpha = \Delta\{D_{g_\alpha(\gamma)}:\gamma\in\omega_1\}$, so every $C_\alpha$ is a club subset of $\aleph_1$.
Let $\mathcal{C}=\{C_\alpha:\alpha\in\mu^+\}$.
Choose a set $A\subseteq\mu^+,|A|=\aleph_1$ such that $\bigcap\{C_\alpha:\alpha\in A\}$ is a club of $\aleph_1$.
Fix an ordinal $\delta\in\omega_1$ such that $|\{g_\alpha(\gamma):\alpha\in A, \gamma<\delta\}|=\aleph_1$.
Define $C=\bigcap\{C_\alpha:\alpha\in A\} - (\delta+1)$, so $C$ is still a club of $\aleph_1$.
Let $D=\bigcap\{D_{g_\alpha(\gamma)}:\alpha\in A,\gamma<\delta\}$, so $D$ is bounded in $\omega_1$.
But if $\zeta\in C$ then for every $\alpha\in A$ and every $\gamma<\delta$ we see that $\zeta\in D_{g_\alpha(\gamma)}$, since $\zeta>\delta>\gamma$.
This means that $C\subseteq D$, a contradiction.

\hfill \qedref{thmcofgp}

For placing the above theorem in a familiar combinatorial land we must analyze the conditions under which the assumption on $A$ can be made.
Let $(f_\alpha:\alpha<\delta)$ be a scale with respect to $J=J^{\rm bd}_\omega$ (this is our typical case).
We shall say that $\gamma$ is a good point if there is an unbounded $u\subseteq\gamma$ and $n\in\omega$ such that $\langle f_\alpha(m):\alpha\in u\rangle$ is strictly increasing for every $m\geq n$.
A scale $(f_\alpha:\alpha<\delta)$ is good iff every $\gamma\in\delta$ of uncountable cofinality is a good point.

\begin{claim}
\label{clmvgs} Assume that $\omega=\cf(\mu)<\mu<2^\omega$, and there exists a good scale $(f_\alpha:\alpha\in\mu^+)$ with respect to some product of regular cardinals and $J=J^{\rm bd}_\omega$. \newline
Then $\mathfrak{gp}\neq\mu$.
\end{claim}

\par\noindent\emph{Proof}. \newline
For every $\beta\in\mu^+$ such that $\cf(\beta)>\omega$ let $u_\beta$ be an unbounded subset of $\beta$ as guaranteed by the assumption that $(f_\alpha:\alpha\in\mu^+)$ is a good scale.
Assume that $A\subseteq\mu^+,|A|=\aleph_1$ and $\beta=\sup(A)$.
Without loss of generality, $\cf(\beta)=\omega_1$.
By shrinking $A$ and $u_\beta$ if needed we may assume that the elements of $A$ and $u_\beta$ are cofinally interleaved.

For every $\alpha\in A$ choose $\delta^\alpha_0,\delta^\alpha_1\in u_\beta$ such that $\delta^\alpha_0<\alpha$ and $\delta^\alpha_1$ is the first element of $u_\beta$ such that $\alpha<\delta^\alpha_1$.
Choose $n=n(\alpha)\in\omega$ such that $f_{\delta^\alpha_0}(m)<f_\alpha(m) <f_{\delta^\alpha_1}(m)$ for every $m\in[n,\omega)$.
By shrinking $A$ once more we may assume that for every $\alpha\in A$ we have $n(\alpha)=n_0$ for some fixed $n_0\in\omega$.

Since $(f_\alpha:\alpha\in\mu^+)$ is a good scale, and $\cf(\beta)>\omega$, there exists $n_1\in\omega$ such that $\langle f_\delta(m):\delta\in u_\beta\rangle$ is strictly increasing for every $m\in[n_1,\omega)$.
Set $n=\max\{n_0,n_1\}$.
It follows that $|\{f_\alpha(\ell):\alpha\in A, \ell\in[n,\omega)\}|=\aleph_1$, as required.

\hfill \qedref{clmvgs}

The above claim enables us to derive conclusions about the cofinality of $\mathfrak{gp}$ from a weak form of the square principle.
It is known that $\square^*_\mu$ implies the existence of a good scale $(f_\alpha:\alpha\in\mu^+)$.
Actually, it gives a bit more. It is shown in \cite{MR2160657} that a better scale can be obtained from the weak square.

\begin{corollary}
\label{corweaks} Assume $\square^*_\mu$ for every $\mu>\cf(\mu)=\omega$ below $2^\omega$. \newline
Then $\cf(\mathfrak{gp})>\omega$.
\end{corollary}

\hfill \qedref{corweaks}

The above conclusions are stated with respect to the possibility that $\cf(\mathfrak{gp})=\omega$, but the real point here seems to be different.
If $\mu>\cf(\mu)=\omega$ and ${\rm pp}(\mu)$ is large then mild assumptions rule out the possibility of $\mathfrak{gp}=\nu$ for many cardinals $\nu\in[\mu,{\rm pp}(\mu))$.
The following is a typical statement:

\begin{theorem}
\label{thmpp} Assume that:
\begin{enumerate}
\item [$(a)$] $\omega=\cf(\mu)<\mu<2^\omega$.
\item [$(b)$] $(\mu_n:n\in\omega)$ is a sequence of regular cardinals, $\mu=\bigcup_{n\in\omega}\mu_n$ and $J\supseteq J^{\rm bd}_\omega$.
\item [$(c)$] $\lambda={\rm tcf}(\prod_{n\in\omega}\mu_n,J)$ and $(f_\alpha:\alpha\in\lambda)$ is a scale which witnesses this fact.
\item [$(d)$] $\mu\leq\nu<\lambda$.
\item [$(e)$] For every $A\in[\lambda]^{\aleph_1}$ it is true that $|\{f_\alpha(n):\alpha\in A, n\in\omega\}|=\aleph_1$.
\end{enumerate}
Then $\mathfrak{gp}\neq\nu$. \newline
\end{theorem}

\par\noindent\emph{Proof}. \newline
Follow the proof of Theorem \ref{thmcofgp}, upon replacing $\mu$ by $\nu$ and $\mu^+$ by $\lambda$.

\hfill \qedref{thmpp}

It follows that the consistency of $\cf(\mathfrak{gp})=\omega$ requires strong assumptions.
We turn to $\mathfrak{gp}_\kappa$ when $\kappa=\cf(\kappa)>\aleph_0$.
If $\kappa>\aleph_0$ then theoretically it may occur that $\cf(\mathfrak{gp}_\kappa)<\kappa$.
Weak assumptions about good scales imply $\cf(\mathfrak{gp}_\kappa)>\kappa$ by a straightforward generalization of the previous statements.
The purpose of the following theorem is to show, in ZFC, that $\cf(\mathfrak{gp}_\kappa)$ cannot drop down dramatically without any special assumption.

\begin{theorem}
\label{thmgpkappa} Martin's maximum and the cofinality of $\mathfrak{gp}_\kappa$. \newline
Assume that $\kappa=\theta^{+3}$. \newline
Then $\cf(\mathfrak{gp}_\kappa)>\theta$, and moreover if $\kappa<\mu$ and $\cf(\mu)\leq\theta$ then $\mathfrak{gp}_\kappa\notin[\mu,{\rm pp}(\mu))$. \newline
Similarly, Martin's maximum implies that $\cf(\mathfrak{gp}_{\theta^+})>\theta$ for every infinite cardinal $\theta$.
\end{theorem}

\par\noindent\emph{Proof}. \newline
Fix any $\mu\in(\kappa,2^\kappa)$ such that $\cf(\mu)\leq\theta$.
For notational simplicity we assume that $\cf(\mu)=\theta$, the proof of the case in which $\cf(\mu)<\theta$ is identical.
Let $\lambda\in[\mu,{\rm pp}(\mu))$ be any cardinal, and assume toward contradiction that $\mathfrak{gp}_\kappa=\lambda$.
Notice that $\lambda^+\leq{\rm pp}(\mu)$ and hence there is a set $\mathfrak{a} = \{\lambda_i:i\in\theta\}\subseteq{\rm Reg}\cap\mu$ such that $\sup(\mathfrak{a})=\mu$ and ${\rm tcf}(\prod\mathfrak{a},J)=\lambda^+$ for some $J\supseteq J^{\rm bd}_\kappa$.

Choose an increasing cofinal sequence $\bar{f} = (f_\alpha:\alpha\in\lambda^+)$ which exemplifies this fact and satisfies $S^{\rm gd}_\theta(\bar{f}) =_{\mathcal{D}_\lambda} S^\lambda_\theta$.
Such a sequence of functions exists by results of Shelah, see e.g. \cite{MR2608402}, footnote 5.
By shrinking $\bar{f}$ if needed we may assume that if $A\subseteq\lambda^+, |A|=\kappa^+$ then $|\{f_\alpha(i):\alpha\in A, i\in\theta\}|=\kappa^+$.
Now the argument of the proof of Theorem \ref{thmcofgp} leads to a contradiction.
The additional part under Martin's maximum follows from \cite{MR2608402}.

\hfill \qedref{thmgpkappa}

Chang's conjecture $(\mu^+,\mu)\twoheadrightarrow(\aleph_1,\aleph_0)$ implies the existence of sets $A\subseteq\mu^+, |A|=\aleph_1$ such that $|\{f_\alpha(n): \alpha\in A,n\in\omega\}|=\aleph_0$ for some scale $(f_\alpha:\alpha\in\mu^+)$ in the relevant context.
But this fact by itself doesn't imply immediately $\cf(\mathfrak{gp})=\omega$.
The point is that even under strong instances of Chang's conjecture we might be able to find a set $A$ of size $\aleph_1$ for which $|\{f_\alpha(n): \alpha\in A,n\in\omega\}|=\aleph_1$.
If we can make sure that such $A$ gives rise to a collection of $\aleph_1$-many clubs with unbounded intersection in $\omega_1$ then $\cf(\mathfrak{gp})>\omega$.
One way to guarantee the above property is captured in the following definition:

\begin{definition}
\label{defstronggp} Strong Galvin's property. \newline
A collection $\{C_\alpha:\alpha\in\lambda\}$ of clubs of $\kappa=\cf(\kappa)>\aleph_0$ exemplifies the strong Galvin's property iff there exists a sequence $\langle H_\varepsilon:\varepsilon<\kappa\rangle$ so that $H_\varepsilon\in[\lambda]^\lambda$ for every $\varepsilon<\kappa$, and for every $h\in\prod_{\varepsilon<\kappa}H_\varepsilon$ the set $\bigcap\{C_{h(\varepsilon)}:\varepsilon<\kappa\}$ is a club subset of $\kappa$.
\end{definition}

The proof of Galvin from \cite{MR0369081} as well as the proof under the weak diamond in \cite{MR3604115} give the strong property for $\lambda=(2^\omega)^+$ and every family of $\lambda$ clubs in $\omega_1$.
Let us prove that this strong property implies $\cf(\mathfrak{gp})>\omega$.

\begin{claim}
\label{clmsgp} Assume that $\mathfrak{gp}_\kappa=\mu$. \newline
If every family $\{C_\alpha:\alpha\in\mu^+\}$ of clubs of $\kappa^+$ has the strong Galvin's property then $\cf(\mu)>\kappa$.
\end{claim}

\par\noindent\emph{Proof}. \newline
Assume toward contradiction that $\theta=\cf(\mu)\leq\kappa$.
Choose a collection $\mathcal{D}=\{D_\alpha:\alpha<\mu\}$ which exemplifies $\mathfrak{gp}_\kappa=\mu$.
Fix a scale $(f_\alpha:\alpha\in\mu^+)$ in $(\prod_{\eta\in\theta}\mu_\eta,J^{\rm bd}_\theta)$.
By chopping an initial segment and thinning out the sequence we can assume that $\kappa<\mu_0$ and $\zeta<\eta\Rightarrow\mu_\zeta<\mu_\eta$.
Denote the ideal $J^{\rm bd}_\theta$ by $J$.

For every $\alpha\in\mu^+$ let $C_\alpha = \bigcap\{D_{f_\alpha(\eta)}:\eta\in\theta\}$, and let $\mathcal{C} = \{C_\alpha:\alpha\in\mu^+\}$.
By our assumptions there is a sequence of sets $\langle H_\varepsilon:\varepsilon<\kappa^+\rangle$, each $H_\varepsilon$ belongs to $[\mu^+]^{\mu^+}$ and for every $h\in\prod_{\varepsilon<\kappa^+}H_\varepsilon$ the set $\bigcap\{C_{h(\varepsilon)}:\varepsilon<\kappa^+\}$ is a club of $\kappa^+$.

We define a set $A\subseteq\mu^+,|A|=\kappa^+$ as follows.
By induction on $\varepsilon<\kappa^+$ we choose an ordinal $\alpha_\varepsilon\in H_\varepsilon$ such that $\sup\{f_{\alpha_\zeta}: \zeta<\varepsilon\}<_J f_{\alpha_\varepsilon}$.
Since $\kappa<\mu_\eta$ for every $\eta$, the function $\sup\{f_{\alpha_\zeta}: \zeta<\varepsilon\}$ is an element of the product.
The existence of $f_{\alpha_\varepsilon}$ follows from the fact that $|H_\varepsilon|=\mu^+$ and hence $\{f_\alpha:\alpha\in H_\varepsilon\}$ is cofinal in the product $(\prod_{\eta\in\theta}\mu_\eta,J)$.

Let $A=\{\alpha_\varepsilon:\varepsilon\in\kappa^+\}$.
By the choice of each $\alpha_\varepsilon$ we see that $|\{f_\alpha(\eta): \alpha\in A, \eta\in\theta\}|=\kappa^+$.
As before, let $C=\bigcap\{C_\alpha:\alpha\in A\}, D=\bigcap\{D_{f_\alpha(\eta)}: \alpha\in A, \eta\in\theta\}$, so $C=D$.
However, $D$ is bounded in $\kappa^+$ while $C$ is a club of $\kappa^+$, a contradiction.

\hfill \qedref{clmsgp}

The above claim invites an additional parameter within the definition of $\mathfrak{gp}$.
For $\kappa=\cf(\kappa)>\aleph_0$ and $\mu\geq\kappa$ we shall say that a family $\{C_\alpha:\alpha\in\lambda\}$ of clubs of $\kappa$ has $\mu$-Galvin's property iff there is a sequence $\langle H_\varepsilon:\varepsilon<\kappa\rangle, H_\varepsilon\in[\lambda]^\mu$ for every $\varepsilon<\kappa$, such that for every $h\in\prod_{\varepsilon<\kappa}H_\varepsilon$ the set $\bigcap\{C_{h(\varepsilon)}: \varepsilon<\kappa\}$ is a club subset of $\kappa$.

\begin{definition}
\label{defgpparameter} $\mathfrak{gp}_\kappa$ with a parameter. \newline
Let $\kappa$ be a regular cardinal, $\mu\in[\kappa^+,(2^\kappa)^+]$. \newline
The characteristic $\mathfrak{gp}_\kappa$ with parameter $\mu$ is the minimal cardinal $\lambda$ such that every collection $\{C_\alpha:\alpha\in\lambda^+\}$ of clubs of $\kappa^+$ has $\mu$-Galvin's property. \newline
If $\mu=(2^\kappa)^+$ then we call it the strong $\mathfrak{gp}_\kappa$ and denote it by $\mathfrak{sgp}_\kappa$.
\end{definition}

In this light, Claim \ref{clmsgp} says that the cofinality of $\mathfrak{sgp}_\kappa$ is always larger than $\kappa$.
Therefore, if we can force $\cf(\mathfrak{gp}_\kappa)=\kappa$ then we can distinguish $\mathfrak{sgp}_\kappa$ from $\mathfrak{gp}_\kappa$.
This can be phrased in a general manner:

\begin{question}
\label{qsgp} Is it consistent that $\mathfrak{gp}$ is strictly smaller than $\mathfrak{sgp}$? \newline
More generally, for $\kappa=\cf(\kappa)$ and $\kappa^+\leq\mu_0<\mu_1\leq(2^\kappa)^+$, is it consistent that $\mathfrak{gp}_\kappa$ with parameter $\mu_0$ is strictly smaller than $\mathfrak{gp}_\kappa$ with parameter $\mu_1$?
\end{question}

In closing this paper we address another problem from \cite{MR3604115}, about the connection between the proper forcing axiom and Galvin's property. Here we refer to the original property proved by Galvin under the continuum hypothesis, namely that every family of $\aleph_2$ many clubs of $\aleph_1$ admits a sub-family of size $\aleph_1$ whose intersection is a club of $\aleph_1$.
We shall prove that consistently this property fails under PFA.
For this end, we shall use the usual way to force PFA from a supercompact cardinal.
It is plausible that the negation of Galvin's property follows from the PFA no matter how it is forced.

\begin{theorem}
\label{thmpfa} PFA and Galvin's property. \newline
It is consistent that PFA holds and there is a family $\{B_\alpha:\alpha\in\omega_2\}$ of clubs of $\omega_1$ for which any subfamily of size $\aleph_1$ has finite intersection.
\end{theorem}

\par\noindent\emph{Proof}. \newline
We denote the ground model ahead of forcing the proper forcing axiom by $V$.
Let $\kappa$ be supercompact in $V$, and force PFA in the usual way described first by Baumgartner.
In the generic extension, $\kappa$ becomes $\aleph_2$ and the PFA holds.
Let $\mathbb{B}$ be Baumgartner's forcing from \cite{MR776640} to add a club of $\aleph_1$ with finite conditions as defined in $V$.
Let $f$ be the Laver function employed in the iteration which forces PFA.
Denote the iteration which forces PFA by $\langle\mathbb{P}_\alpha,\name{\mathbb{Q}_\beta}:\alpha\leq\kappa,\beta<\kappa\rangle$, and fix a generic subset $G\subseteq\mathbb{P}_\kappa$.

Let $A=\{\zeta\in\kappa:\ \Vdash_{\mathbb{P}_\zeta} f(\zeta)=\check{\mathbb{B}}\}$.
If $\zeta\in A$ then $\mathbb{Q}_\zeta$ adds a club of $\aleph_1$, denoted by $\name{\mathbb{B}}_\zeta$.
We claim that the sequence $\langle\name{\mathbb{B}}_\zeta:\zeta\in A\rangle$ exemplifies the negation of Galvin's property in the generic extension by $\mathbb{P}_\kappa$.
Remark that in the generic extension this is a set of size $\aleph_2$, so this will accomplish the proof.

Assume, therefore, that $\name{X}$ is a name of an element of $[A]^{\aleph_1}$.
Let $\name{\tau} = \bigcap_{\zeta\in\name{X}}\name{\mathbb{B}}_\zeta$.
Assume toward contradiction that $\tau$ is forced to be uncountable.
Fix an element $p$ of $\mathbb{P}_\kappa$, with the goal to extend $p$ to a condition which forces the opposite statement.
Since $p$ is arbitrary, we will be done.
For a sufficiently large regular $\chi$ we fix $M\prec\mathcal{H}(\chi)$ such that $p,\name{X},\mathbb{P}_\kappa,\name{\tau}\in M$ and $M$ is countable.
By properness, let $q$ be an $(M,\mathbb{P}_\kappa)$-generic condition such that $p\leq q$.
Define:
$$
\name{\sigma} = \{(r,\gamma)\in M: r\Vdash\check{\gamma}\in\name{\tau}\}.
$$
The focal point is that $q\Vdash\name{\sigma}=\name{\tau}\cap M$.
We prove the inclusion of $\name{\tau}\cap M$ in $\name{\sigma}$, the opposite inclusion being proved exactly in the same manner.
Assume, therefore, that $\gamma\in M$.
For the current direction we assume that there is a condition $s\geq q$ such that $s\Vdash\check{\gamma}\in\name{\tau}$.
Since $q$ is $(M,\mathbb{P}_\kappa)$-generic and the set of conditions which decide the statement $\check{\gamma}\in\name{\tau}$ is dense, one can find $t\in M$ such that $t\parallel s \wedge t\Vdash\check{\gamma}\in\name{\tau}$.
By the definition of $\name{\sigma}$ we see that $t\Vdash\check{\gamma}\in\name{\sigma}$ and hence some extension of $s$ forces this statement as well, giving the desired inclusion.
The opposite inclusion is proved with respect to the opposite statement, namely $s\Vdash\check{\gamma}\notin\name{\tau}$.
Remark that $q\Vdash\name{\sigma}=\name{\tau}\cap M$ along with the assumption toward contradiction imply that $\name{\sigma}$ is forced by $q$ to be an infinite set.

Let $\eta=\sup(\name{X})$.
We may assume that this fact is forced by $q$.
By cutting off the upper part of $X$ if necessary, we may assume that $\cf^V(\eta)>\omega$.
Let $\delta=\sup(M\cap\eta)$.
Recall that $M$ is countable, and hence $\delta<\eta$.
Choose a condition $q^+\geq q$ such that $q^+\Vdash\check{\gamma}\in\name{X}\wedge\gamma\in[\delta,\eta)$.
We shall argue that $q^+\Vdash \neg(\name{\sigma}\subseteq \name{\mathbb{B}}_\gamma\cap M)$, thus arriving at a contradiction.

Our strategy is to prove the above statement in a larger universe, and then using absoluteness in order to show that this holds already in the generic extension by $\mathbb{P}_\kappa$.
Let $\mathbb{T}$ be the termspace forcing of $\mathbb{P}_{[\eta,\kappa)}$ with respect to $\mathbb{P}_{[\delta,\eta)}$ as computed in the universe obtained by the generic extension with $\mathbb{P}_\delta$.
Let $\mathbb{R} = \mathbb{P}_\delta\ast(\mathbb{P}_{[\delta,\eta)}\times\mathbb{T})$.
Observe that $\mathbb{R}$ has a natural projection onto $\mathbb{P}_\kappa$, namely $\pi(p,q,r)=p^\frown q^\frown r$.
In particular, there are suitable conditions $(a,b,c)$ such that $\pi(a,b,c)=q^+$.
Choose a generic set $H$ for the forcing notion $\mathbb{P}_\delta\ast\mathbb{T}$, such that $(a,c)\in H$.
Notice that the name $\name{\sigma}$ is interpreted by $H$ since there are no coordinates of $\name{\sigma}$ in the interval of $\mathbb{P}_{[\delta,\eta)}$ by the definition of $\name{\sigma}$.

We claim that in $V[H]$ it is true that $\varnothing_{\mathbb{P}_{[\delta,\eta)}} \Vdash \neg(\name{\sigma}\subseteq \name{B}_\gamma\cap M)$.
Remark that maybe $\aleph_1$ is collapsed in $V[H]$ (though not in $V[G]$), but our claim is still valid.
Indeed, if $p$ is a finite condition in the forcing which adds $\mathbb{B}_\gamma$ then we can extend $p$ by adding a pair which will be forced to be in $\name{\sigma}$ but not in $\name{\mathbb{B}_\gamma}\cap M$.
For this end, let $\alpha=\sup(\name{\sigma})$.
Since $p$ is finite, there is a maximal pair $(\nu,\xi)$ in $p$.
Choose an ordinal $\rho\in\name{\sigma}\cap(\xi,\alpha)$ such that $\rho+1$ will not be in $\name{\sigma}$ (such an ordinal always exists) and add the pair $(\nu+1,\rho+1)$ to $p$.
This shows that $\neg(\name{\sigma}\subseteq \name{\mathbb{B}}_\gamma\cap M)$ in the generic extension by $H$.
Since this statement is $\Delta_0$ we conclude that $V[G]\models\neg(\name{\sigma}\subseteq \name{\mathbb{B}}_\gamma\cap M)$ as well, so we are done.

\hfill \qedref{thmpfa}

We make the comment that an additional argument shows that Galvin's property fails in a strong sense in the above model, namely the intersection of any family of $\aleph_1$ Baumgartner clubs is finite.
The important feature of the usual iteration to force PFA that we used is the fact that the set $A$ of names for Baumgartner clubs along the iteration is unbounded in $\kappa$.

\newpage

\bibliographystyle{amsplain}
\bibliography{arlist}

\end{document}